\magnification 1200
\def\R{{\rm I\kern-0.2em R\kern0.2em \kern-0.2em}}
\def\N{{\rm I\kern-0.2em N\kern0.2em \kern-0.2em}}
\def\P{{\rm I\kern-0.2em P\kern0.2em \kern-0.2em}}
\def\B{{\rm I\kern-0.2em B\kern0.2em \kern-0.2em}}
\def\C{{\bf \rm C}\kern-.4em {\vrule height1.4ex width.08em depth-.04ex}\;}

\def\D{{\Delta}}

\def\z{{\zeta}}
\def\cC{{\cal C}}
\def\cP{{\cal P}}

\def\cW{{\cal W}}
\def\cV{{\cal V}}

\font\ninerm=cmr8
\centerline {\bf ANALYTICITY ON FAMILIES OF CIRCLES}
\vskip 4mm
\centerline {Josip Globevnik}
\vskip 4mm
{\noindent \ninerm ABSTRACT\  It is known that if f is a 
continuous function on the
complex plane which extends holomorphically from each circle 
surrounding the origin then f is 
not necessarily holomorphic. In the paper we prove that if, 
in addition, f extends holomorphically 
from each circle belonging to an open family of circles which 
do not surround the origin then f 
is holomorphic.} 
\vskip 4mm
\bf 1.\ Introduction and the main result \rm 
\vskip 2mm
Write $\D (a,\rho)= \{\z\in\C\colon\ |\z -a|<\rho\}$ 
and $\D = \D (0,1)$.  If 
$0<r_1<r_2<\infty$ write $A(a,r_1,r_2)=\{\z\in\C\colon\ r_1\leq |\z-a|\leq r_2\}$. 
We say 
that a continuous function on $b\D (a,\rho)$ extends holomorphically 
from $b\D (a,\rho)$ 
if it has a continuous extension to $\overline\D (a,\rho)$ which 
is holomorphic on $\D (a,\rho)$.

A family $\cC$ of circles is called a \it test family for holomorphy \rm (on $\C $ ) 
if every 
continuous function on $\C $ that extends holomorphically from each circle in
$\cC$ is holomorphic 
on $\C $. We will consider open families of circles, that is, families of 
the form 
$\{ b\D (a,\rho)\colon\ (a,\rho)\in\cP\}$ where $\cP$ is an open subset of 
$\C\times (0,\infty)$. 

There are large families of circles that are not test families for holomorphy.
For instance, 
the function
$$
f(z)= 
\left\{\eqalign{&z^2/\overline z\ \ (z\in\C\setminus\{ 0\})\cr
&0\ \ \ \ \ \ \ (z=0)\cr}\right.
$$
is continuous on $\C$ and extends holomorphically from each circle
that surrounds the origin,
yet $f$ is not holomorphic. This shows that the family of all circles 
that surround the origin is
not a test family for holomorphy. In the present paper we prove
 that the family of all circles that surround the origin is a maximal 
open family that is not a test family for holomorphy:
\vskip 2mm
\noindent\bf Theorem 1.1\ \ \it Let $f$ be a continuous function on $\C\setminus \{ 0\}$ 
which extends holomorphically from each circle that surrounds the origin.
Suppose that, in addition, 
$f$ extends holomorphically from each circle belonging to a nonempty open 
family of circles 
that do not surround the origin. Then $f$ is an entire function, that is, $f$
is a holomorphic function 
on $\C\setminus \{ 0 \}$ which has a removable singularity at $0$. \rm 
\vskip 2mm
\noindent We prove Theorem 1.1 in the first part of the paper. In the second 
part we look at special cases  of nonholomorphic continuous functions 
on $\C\setminus\{ 0\}$ which 
extend
holomorphically from every circle surrounding the origin. In particular, 
we consider functions constant on lines passing 
through the origin and functions constant on rays passing through the origin. 

To prove Theorem 1.1 we use a new approach to the holomorphic 
extension problem for circles which was introduced in [AG] and further developed 
in [G3]. We describe this new approach. In [AG] we studied rational 
functions of two real variables \ $f(z) = P(z,\overline z)/Q(z,\overline z)$ where 
$P, Q$ are polynomials. We noticed that $f|b\D (a,\rho)$ has a unique meromorphic extension 
to $\D (a,\rho )$ given by 
$$
f^\ast (z) = {{P(z,\overline a+\rho^2/(z-a))}\over{Q(z,\overline a+\rho^2/(z-a))}}
\eqno (1.1)
$$
so to say that $f$ extends holomorphically from $b\D (a,\rho)$ means that 
$f^\ast $ has no singularities in $\D (a,\rho )$. Given $a\in\C$ and $\rho >0$ 
we introduced
$$
\Lambda _{a,\rho} = \{ (z,w)\in\C ^2\colon (z-a)(w-\overline a)=\rho^2, 0<|z-a|<\rho\} 
$$
a closed complex submanifold of $\C^2\setminus\Sigma$,\ 
attached to the real two-plane $\Sigma =\{ (z, \overline z)\colon z\in\C\}$  along the circle 
$b\Lambda _{a,\rho} = \{(z,\overline z)\colon z\in b\D(a,\rho)\}$. The holomorphic extension 
of $f$ from $b\D (a,\rho)$ to $\D (a,\rho)$ is then given as the restriction of the rational 
function $P(z,w)/Q(z,w)$ of two complex variables to the complex manifold 
$\Lambda _{a,\rho}$ as seen from (1.1). 

In [G3] the varieties $\Lambda _{a,\rho}$ were used to formulate the holomorphic extension 
problem for general continous functions as a problem in $\C^2$, starting from the trivial
observation that 
a continuous function $f$ on $b\D (a,\rho)$ extends holomorphically to $\D (a,\rho)$ if 
and only if the function $F(z,\overline z)= f(z)$ defined on $b\Lambda _{a,\rho}$ has a 
bounded 
continuous extension to $\Lambda _{a,\rho}\cup b\Lambda _{a,\rho}$ which is 
holomorphic on $\Lambda _{a,\rho}$. If $A= A(b,r_1,r_2)$ it
was shown that 
the union $\Omega (A)$ of all $\Lambda _{a,\rho}$ such that $b\D (a,\rho)\subset \hbox{Int}A$ 
surrounds $b$ is a 
wedge domain attached to $\Sigma $ along $\tilde A = \{ (z,\overline z)\colon z\in A\}$. If 
$f$ is a continuous function on $ A$ which extends holomorphically from each 
circle $b\D (a,\rho)\subset A$ surrounding $b$, then by an old result of the author [G1], the function 
$f$ is 
the uniform limit of a sequence of polynomials in $z-b$ and $1/(\overline z-\overline b)$ 
which implies that the function $F(z,\overline z)=f(z)\ (z\in A)$ has a bounded 
continuous extension to $\Omega (A)\cup b\Omega (A)$ which is holomorphic on 
$\Omega (A)$. So, roughly 
speaking, continuous functions $f$ extendible holomorphically from open families of circles 
are the functions of the form $F(z,\overline z)$ which are the boundary values 
of bounded holomorphic 
functions $F$ on wedge domains attached to $\Sigma $. Clearly $f$ is holomorphic 
if and only if $F$ 
depends only on the first variable. At this point one can apply standard tools 
of several 
complex variables. We do this to prove Theorem 1.1. 

A different formulation of the holomorphic 
extension problem as a problem in $\C^2$ has been used by A.\ Tumanov [T] for 
continuous functions $f$ on the strip $\{ z\in\C\colon |\Im z|\leq 1\}$ which extend 
holomorphically from each circle $b\D (t,1),\ t\in\R$. Tumanov defines 
a function $F$ on 
$M=\{ (\z +t,\z )\colon\ \z\in\overline\D,\ t\in\R\}$, the disjoint 
union of translates of the
disc $\{(\z , \z\}\colon\ \z\in\D\}$ in such a 
way that for each $t\in\R$,\ the function 
$\z\mapsto F(\z +t,\z)\ (\z\in\overline\D)$ is the continous extension of $\z\mapsto 
f(\z +t)$ to $\overline\D$ which is holomorphic on $\D $. In particular, $F(\z +t,\z ) =
f(\z + t)\ (\z\in b\D , t\in \R)$. He observes that 
since $F(z,w)=F(z,-1/w)\ (w\in b\D )$ one can extend $F$ to 
$\tilde M = \{ (\z + t, -1/\z)\colon\  
t\in\R, \z\in \overline \D\setminus \{ 0\}\}$ by $F(z,w)\equiv F(z,-1/w)$ to get a 
continuous CR function on the CR manifold $M\cup\tilde M$. He then constructs
analytic discs attached to $M\cup\tilde M$ and uses 
the Baouendi-Treves approximation theorem, the 
edge of the wedge theorem and the continuity principle to prove that $F$ does not depend on 
the second variable, that is, that $f$ is holomorphic on $\{\z\in\C\colon\ |\Im \z |<1\}$. 
\vskip 4mm
\bf 2.\ Varieties $\Lambda _{a,\rho}$ and domains $\Omega (A)$ \rm
\vskip 2mm
Let $0<r_1<r_2<\infty$ and let $a\in\C$. Denote by $\Omega (A(a,r_1,r_2))$ 
the union of 
all $\Lambda _{b,\rho}$ such that $b\D (b,\rho)\subset\hbox{Int}A(a,r_1,r_2)$
surrounds 
$a$. The set $\Omega (A(a,r_1,r_2))$ is an unbounded open connected 
subset of $\C^2\setminus\Sigma $ 
which is attached to $\Sigma $ along 
$\tilde A(a,r_1,r_2)=\{(z,\overline z)\colon\ z\in A(a,r_1,r_2)\}$. We shall 
need the following
\vskip 2mm
\noindent\bf Theorem 2.1\ [G3]\it \ \ Let $f$ be a continuous
function on $A(a,r_1,r_2)$. 
The following are equivalent

(i)\ $f$ extends holomorphically from each circle $b\D (b,\rho)\subset A(a,r_1,r_2)$ 
which surrounds the point $a$

(ii)\ the function $F(z,\overline z)=f(z)$ defined on 
$\{(z,\overline z)\colon\ z\in A(a,r_1,r_2)\}$ 
extends to a bounded continuous function on 
$\Omega (A(a,r_1,r_2))\cup b\Omega (A(a,r_1,r_2))$ which is holomorphic on $\Omega 
(A(a,r_1,r_2))$. \rm 
\vskip 2mm
We list some simple properties of $\Lambda _{a,\rho}$ and $\Omega (A)$.
The proofs are 
elementary. They can be found in [G3]. The proof of Proposition 2.1 can be found also in the 
earlier paper [AG]. 
\vskip 2mm
\noindent\bf Proposition 2.1\ \ \it Let $(z,w)\in\C^2\setminus\Sigma$. 
Then $(z,w)\in\Lambda _{
a,R}$ if and only if there is a $t>0$ such that 
$a=z+t(z-\overline w)$ and 
$R=\sqrt{t(t+1)}|z-\overline w|$. In fact, given $R>0$ we have
$$
a=z+2^{-1}[\sqrt{1+4R^2/|z-\overline w|^2}-1](z-\overline w) .
\eqno (2.1)
$$
\vskip 1mm\rm 
Note that the two dimensional subspace perpendicular to the Lagrangian 
two-plane $\Sigma$ is $i\Sigma = \{(z,-\overline z)\colon\ z\in\C\}$. 
Our next lemma tells how 
a variety $\Lambda _{a,\rho}$ intersects the two dimensional planes
perpendicular to $\Sigma $: 
\vskip 2mm
\noindent \bf Proposition 2.2\ \ \it Let $z\in\C,\ t>0$ and $\varphi \in\R$. Then 
$(z,\overline z)+(te^{i\varphi},-te^{-i\varphi})\in\Lambda _{a,R}$ if and only if 
$a=z+\sqrt{t^2+R^2}e^{i\varphi}$. \rm 
\vskip 2mm
\noindent\bf Proposition 2.3\ \ \it Let $A=A(a,r_1,r_2)$. 
Then $\Omega (A)$ 
is an unbounded open connected subset of $\C^2\setminus \Sigma$ attached to $\Sigma$
 along $\{ (z,\overline z)\colon\ z\in\hbox{Int}A\}$. If $\gamma = (r_1+r_2)/2$ then 
 $b\Omega (A)$ consists of $\tilde A = \{(z,\overline z)\colon\ z\in A\}$ 
 together with all 
 $\Lambda _{b,\rho}$ associated with those $b\D (b,\gamma )\subset A$ 
 which are tangent to 
 both $b\D (a,r_1)$ and $b\D (a,r_2)$. Further, $\Omega (A)$ is a 
 disjoint union 
 of $\Lambda _{b,\gamma} $ such that $b\D (b,\gamma)\subset\hbox{Int}A$. \rm
\vskip 2mm

Let $z_0\in\hbox{Int}A$. Let $\Gamma (z_0)\subset \hbox{Int}A$ be the
circle of radius $\gamma =
(r_1+r_2)/2$. which passes through $z_0$ and whose center $b(z_0)$ lies 
on the line through $a$ and $z_0$. For 
each $\varphi \in\R$ define $T_\varphi (z) = 
z_0+e^{i\varphi}(z-z_0)$. $T_\varphi $ is the rotation with center 
$z_0$ for the angle $\varphi$. There is a $\delta (z_0),\ 0<\delta (z_0)<\pi /2$, such 
that $T_\varphi (\Gamma (z_0))\subset\hbox{Int}A\ \ (-\delta (z_0)<\varphi <\delta (z_0))$ 
and such that both $T_{\delta (z_0)}(\Gamma (z_0))$ and $T_{-\delta (z_0)}(\Gamma (z_0))$ meet 
$bA$\ (in fact, they meet both circles that bound $A$). Fix $\varphi ,\ 
-\delta (z_0)<\varphi <\delta (z_0)$. There is 
a $\tau (z_0,\varphi)>0$ such that
$$
T_\varphi (\Gamma (z_0)) + t{{b(z_0)-z_0}\over{|b(z_0)-z_0|}}e^{i\varphi}
\subset\hbox{Int}A\ \ (0\leq t<\tau (z_0,\varphi ))
$$
while 
$T_\varphi (\Gamma (z_0)) + \tau (z_0,\varphi)((b(z_0)-z_0)/
|b(z_0)-z_0|)e^{i\varphi}$ 
meets $bA$ (in fact, it meets both circles 
that bound $A$). For each  $\varphi ,\ 
-\delta (z_0)<\varphi <\delta (z_0)$,\ let $\eta (z_0,\varphi)= 
\sqrt{\tau (z_0,\varphi)^2+2\tau (z_0,\varphi )\gamma }$. 
By Proposition 2.2 we have
$$
(z_0,\overline{z_0})+
\biggl(t{{b(z_0)-z_0}\over{|b(z_0)-z_0|}}e^{i\varphi},
-t{{\overline{b(z_0)}-\overline{z_0}}\over{|b(z_0)-z_0|}}e^{-i\varphi}\biggr)
\in \Omega (A)
$$
provided that $0<t<\eta (z_0,\varphi)$.  Let
$$
D(z_0)=\{ t{{b(z_0)-z_0}\over{|b(z_0)-z_0|}}e^{i\varphi}\colon 
\ 0<t<\eta (z_0,\varphi ),\ 
-\delta (z_0)<\varphi <\delta (z_0)\}.
$$
It is easy to see that the function $\delta $ is continuous on $\hbox{Int}A$ 
and that $\eta $ is a continuous function of $z_0$ and $\varphi $ where it
is defined. For each 
$z_0\in\hbox{Int}A$ the set
$$
(z_0,\overline{z_0})+\{ (\z ,-\overline\z )\colon\ \z\in D(z_0)\}
$$
is contained in $\Omega (A)$.  This proves
\vskip 2mm
\noindent\bf Proposition 2.4\ \ \it Let $z_0\in\hbox{\rm Int}A$. 
There are a neighbourhood $U\subset\Sigma$ 
of $(z_0,\overline{z_0})$, an open convex cone $K\subset \C$ 
with vertex at the origin, containing $\{ t(a-z_0)\colon\ t>0\}$, 
and an $r>0$ such that if 
$$
P = \{ (\z,-\overline \z )\colon\ \z\in K,\ |\z |<r\}
$$
then  $U+P\subset \Omega (A)$.
\vskip 2mm\rm
Fix $R>0$. From (2.1) we get that if $(z,w)\in\Lambda _{a,R}$ then 
$|a|\leq |z|+2R^2/|z-\overline w|$ which, by Proposition 
2.3, implies that given $r_1, r_2,\ 0<r_1<r_2<\infty $
$$
\eqalign{
&\hbox{there are\ }\delta >0\hbox{\ and\ }M<\infty\hbox{\ such that\ }\cr
&\{(z,w)\colon\ |z|\leq \delta,\ |w|\geq M\}\subset\Omega (A(0,r_1,r_2)).\cr}
\eqno (2.2)
$$
\vskip 3mm
\bf 3.\ Functions that extend holomorphically from every 
circle which surrounds the origin \rm 
\vskip 2mm
Suppose that $f$ is a continuous function on $\C\setminus \{ 0\}$ 
which extends holomorphically from each circle that surrounds the origin. 
Define $F$ on 
$\Sigma\setminus \{ (0,0)\} $ by
$$
F(z,\overline z) = f(z).
$$
Then for each $a,\ \rho $ such that $b\D (a,\rho )$ surrounds the origin, the function 
$F|b\Lambda_{a,\rho}$ has a bounded continuous extension to $\Lambda _{a,\rho}\cup 
b\Lambda _{a,\rho}$ which is holomorphic on $\Lambda _{a,\rho}$. 
By Theorem 2.1 we know that 
this defines a holomorphic function $F$ on $\Omega$, the union of all 
$\Lambda _{a,\rho}$ such that $b\D (a,\rho )$ surrounds the
origin. By Theorem 2.1 for each 
$r_1, r_2,\ 0<r_1<r_2<\infty $, the restriction of $F$ to 
$\Omega (A(0,r_1,r_2))$ is bounded and has 
a (bounded) continuous extension to 
$\Omega (A(0,r_1,r_2))\cup b\Omega (A(0,r_1,r_2))$ which, 
on $\tilde A(0,r_1,r_2) = \{ (z,\overline z)\colon\ \z\in A(0,r_1,r_2)\}$ coincides 
with $F(z,\overline z)$.

We now show that 
$$
\Omega = \{(z,w)\colon |w|>|z|\}.
$$
One way to see this is 
by using Proposition 2.2. 
We show this by using Proposition 2.1.

Suppose that $(z,w)\in\C^2\setminus \Sigma$, that is, $w\not= \overline z$.
By Proposition 2.1 
we have $(z,w)\in\Lambda _{a,\rho}$ if and only if 
$$
a= z+t(z-\overline w),\ \ \ \rho = \sqrt{t(t+1)}|z-\overline w|
\eqno (3.1)
$$ for some $t>0$. Let $L$ be the line through $(z+\overline w)/2$
which is perpendicular to 
the line through $z$ and $\overline w$ and let $\Pi\subset\C$ be the open halfplane 
bounded by $L$ which contains 
$z$. It is easy to see that $\Pi$ is the union of all 
$\Delta (a,\rho )$ such that 
$a$ and $\rho $ satisfy (3.1) for some $t>0$. It 
follows that $(z.w)\in\Lambda _{a,\rho}$ 
for some $b\D (a,\rho)$ that surrounds the origin 
if and only if $0\in\Pi$, that is, if and only if $|w|>|z|$.

For each $z\not= 0$ we describe $[(z,\overline z)+i\Sigma ]\cap \Omega $. Recall that 
$i\Sigma = \{(\z,-\overline \z)\colon\ \z\in\C\}$. Write $z=|z|e^{i\alpha}$. Then
 $(z,\overline z)+(\z,-\overline \z)\in\Omega $ if and only 
 if $|z+\z|<|\overline z-\overline\z|$, that is, 
 if and only if $\hbox{Re}(\overline z\z)<0$. 
This happens if and only if $\hbox{Re}(e^{-i\alpha}\z)<0$, that is, if and only if 
 $\z\in e^{i\alpha}\{z\colon\ \hbox{Re}z<0\}$. Let  $L(z)$ be the line through the 
 origin which is perpendicular to the line through $0$ and $z$ and let $P(z)$ be the 
 halfplane bounded by $L(z)$ which does not contain $z$. Then
 $$
 [(z,\overline z)+i\Sigma ]\cap \Omega = (z,\overline z)+
 \{ (\z ,-\overline\z )\colon\ \z\in P(z)\}.
 $$
 Obviously $i\Sigma\cap\Omega=\emptyset$. Thus, $\Omega $ 
 can be written as a disjoint union of 
 halpfplanes 
 $$
 \Omega = \cup_{z\in\C\setminus\{ 0\}}[(z,\overline z)+
 \{ (\z ,-\overline\z )\colon\ \z\in P(z)\}].
 $$
 Further,
 $$
 b\Omega = (i\Sigma)\cup(\cup_{z\in\C\setminus\{ 0\}}
 \{ (\z ,-\overline\z )\colon\ \z\in L(z)\}).
$$

Note that we cannot conclude in general that the function 
$F$ extends continuously to 
$\Sigma \setminus\{(0,0)\}$. If 
$(z,\overline z)\in\Sigma \setminus \{(0,0)\}$ and 
if $(z_n,w_n)\in\Omega,\ 
(z_n,w_n)\rightarrow (z,\overline z)$ then
$\lim_{n\rightarrow\infty}F(z_n,$ $w_n)$ $=F(z,$ $
\overline z)= f(z)$ provided that there are
$r_1, r_2, 0<r_1<r_2<\infty $ such that 
$(z_n,w_n)\in\Omega (A(0,r_1,r_2))$ for all $n$. However, we have the following
\vskip 2mm
\noindent\bf Proposition 3.1\ \ \it Let $f$ and $F$ be as above. Suppose that 
$F(z,\overline z)=f(z)$ has a holomorphic
extension $\Phi $ into an open ball
$B\subset\C^2\setminus\{(0,0)\}$ centered at 
$(z_0,\overline{z_0})\in\Sigma\setminus\{(0,0)\}$.   
Then $\Phi \equiv F$ on 
$B\cap\Omega $. \rm
\vskip 1mm
\noindent\bf Proof.\ \rm  By Proposition 2.4 there are
a neighbourhood $U\subset\Sigma$
 of $(z_0,\overline{z_0})$, an open convex cone 
 $K\subset i\Sigma $ with vertex at the origin 
 and an $\eta >0$ such that if 
 $K_\eta =\{ w\in K,\ |w|<\eta \}$ then $U+K_\eta \subset 
 \Omega (A(0,r_1,r_2)\cap B$ for some $r_1, r_2,  0<r_1<r_2<\infty$  and 
 hence $F(z,\overline z)$ 
 has a continuous extension from $U$ to $U\cup[U+K_\eta ]$ 
 which is holomorphic on $U+K_\eta $. 
 However, such extension is unique and since $\Phi |[U\cup[U+K_\eta]]$
 is such an extension 
 we must have $\Phi\equiv F$ on $U+K_\eta $. Since  $(U+K_\eta )\cap B$ 
 is an open subset of 
 $\Omega\cap B$ and since $\Omega\cap B$ is connected 
 it follows that $\Phi \equiv F$ on 
 $\Omega\cap B$. This completes the proof.
 \vskip 4mm
 \bf 4.\ Intersecting varieties $V_{a,\rho}$ with $\Omega $\rm 
 \vskip 2mm
 Given $a\in\C$ and $\rho >0$ let 
 $$
 V_{a,\rho} = \{(z,w)\colon\ (z-a)(w-\overline a)=\rho^2\}.
 $$
 Thus, $\Lambda_{a,\rho}=\{ (z,w)\in V_{a,\rho}:\ 0<|z-a|<\rho\}$ is one of the 
 two components 
 of $V_{a,\rho}\setminus b\Lambda _{a,\rho}$. 
 We compute $V_{a,\rho}\cap b\Omega $. 
   The equation of $V_{a,\rho}$ is 
 $w=\overline a+\rho^2/(z-a)$, so we compute the 
 intersection of $V_{a,\rho}$ with 
 $b\Omega =\{ (z,w)\colon\ |z|=|w|\}$ by solving 
 $|\overline a+\rho^2/(z-a)|=|z|$. We 
 get  
 $|\overline a(z-a)+\rho^2|=|z|.|z-a|$ so
 $$
 [\overline a(z-a)+\rho^2].[a(\overline z-\overline a)+\rho^2 ]-\rho^2z\overline z=
 z\overline z=
 z\overline z[(z-a)(\overline z-\overline a)-\rho^2].
 $$
 The left hand side equals 
 $$
 \eqalign{
 &a\overline a[(z-a)+\rho^2/{\overline a}].[(\overline z-\overline a)+\rho^2/a]-
 \rho^2[(z-a)(\overline z-\overline a)+a\overline z+\overline a z-a\overline a] =\cr
 &=a\overline a(z-a)(\overline z-\overline a)+\rho^4-\rho^2(z-a)(\overline z-\overline a)+
 a\overline a\rho^2 \cr
 &=(a\overline a-\rho^2)[(z-a)(\overline z-\overline a)-\rho^2]\cr}
 $$
 and the equation becomes 
 $$
 [z\overline z- (a\overline a-\rho^2].[(z-a)(\overline z-\overline a)-\rho^2]= 0.
 $$

If the circle $b\D (a,\rho)$ surrounds the origin, that is, if $|a|<\rho$ then
the set of solutions is 
$b\D (a,\rho)$. The case of interest to us will be the case when $|a|>\rho$, that is, 
when $b\D (a,\rho )$ does not surround the origin. 
In this case the set of solutions is 
$b\D (a,\rho)\cup b\D (0,\sqrt{|a|^2-\rho^2})$.  Note that these two circles 
intersect at right angle. 

Since  
$$
|z-a|^2[|\overline a+\rho^2/(z-a)|^2-|z|^2]
=[(a\overline a-\rho^2 -z\overline z].[(z-a)(\overline z-\overline a)-\rho^2]
\eqno (4.1)
$$
 the point $(z,\overline a+\rho ^2/(z-a))$ belongs to $\Omega $ 
if and only if the 
expression on the left in (4.1) is positive, that is, 
if and only if the expression on the right in (4.1) is 
positive. This happens if and only if either 
$z\in\D (0,\sqrt{|a|^2-r^2})\setminus\overline\D (a,r)$ or 
$z\in\D (a,r)\setminus \overline\D(0,\sqrt{|a|^2-r^2})$. Thus, if
$$
\eqalign{
&D_1(a,r) = \D (a,r)\cap\D(0,\sqrt{|a|^2-r^2})\cr
&D_2(a,r)= \C\setminus\bigl[ \overline\D (a,r)\cup\overline \D(0,\sqrt{|a|^2-r^2})\bigr] \cr
&D_3(a,r)=\D(0,\sqrt{|a|^2-r^2})\setminus \overline\D (a,r) \cr
&D_4(a,r)= \D (a,r)\setminus \overline\D(0,\sqrt{|a|^2-r^2})\cr}
$$
and
$$
\eqalign{&V_i(a,r) = \{ (z,\overline a+\rho ^2/(z-a))\colon\ z\in D_i (a,r)\} \ (1\leq i\leq 3) \cr
&V_4(A,r)=\{ (z,\overline a+\rho ^2/(z-a))\colon\ z\in D_4(a,r)\setminus \{ a\}\} \cr}
$$
then $V_1(a,r)$ and $V_2(a,r)$ are the components of $V_{a,r}\setminus\overline\Omega$ and 
$V_3(a,r)$ and $V_4(a,r)$ are the components of $V_{a,r}\cap\Omega $.  
\vskip 3mm

\bf 5.\ Outline of the proof of Theorem 1.1
\vskip 2mm
\rm
We start with a continuous function $f$ on $\C\setminus\{ 0\}$ which extends 
holomorphically from every circle which surrounds the origin and the associated 
function $F$, holomorphic on $\Omega =\{ (z,w)\colon\ |w|>|z|\}$. Suppose that 
$f$ extends holomorphically from a circle $b\D (b,r)$ that does not surround 
the origin and from all
nearby circles $b\D (a,r)$ with $a$ close to $b$. Then 
$F|b\Lambda _{b,r}$ has a bounded continuous extension $F_1$ 
to $\Lambda _{b,r}\cup b\Lambda _{b,r}$ 
which is holomorphic on $\Lambda _{b,r}$ In particular, $F_1$ 
on $\Lambda _{b,r}\cap \Omega = V_4(b,r)$. We use the edge of the wedge theorem 
as in [G3] to show that, 
since
$f$ extends holomorphically from all nearby circles $b\D (a,\rho)$, on $V_4(b,r)$, 
the function $F_1|V_4(b,r)$ coincides with 
$F|V_4(b,r)$. Since $F_1$ is holomorphic on $\Lambda _{a,r}$, it follows that
$F$ extends holomorphically along $V_{b,r}$ to  $V_1(b,r)$. 
Further, using again the fact that
$f$ extends holomorphically from each circle $b\D (a,r)$ 
where $a$ runs through a neighbourhood 
of $b$ and  repeating the process above with $b$ replaced by $a$ we see that $F$ extends 
holomorphically 
into a neighbourhood $P$ of $\overline{V_1(b,r)}$ in $\C^2$. Now we can apply the 
continuity principle. 
$V_1(b,r)$ is a holomorphically embedded disc which can be continuously 
deformed through a 
family of holomorphically embedded discs into a holomorphically embedded 
disc lying on the $w$-axis which 
contains the origin in its interior, in such a way that
boundaries of all these discs are
contained in $\Omega \cup P$. This implies that $F$ extends holomorphically 
into a neighbourhood of the 
origin. This shows that $f$ can be defined at $0$ so that it becomes 
a real analytic function in a neighbourhood of 
the origin. By a result from [G2] it follows that $f$ is holomorphic 
on $\C$ which will complete the proof.
\vskip 3mm
\bf 6.\ Proof of Theorem 1.1 \rm
\vskip 2mm
Suppose that $f$ extends holomorphically from each circle that surrounds the origin.
We know that there is a holomorphic function 
$F$ on $\Omega $ which, for each $R_1, R_2,\ 0<R_1<R_2$,  has 
a bounded continuous extension to
$\Omega (A(0,R_1,R_2))\cup b\Omega (A(0,R_1,R_2)$ which coincides with 
$F(z,\overline z)=f(z)$ on $\tilde A(0,R_1,R_2)=\{(z,\overline z)\colon\ z\in A(0,R_1,R_2)\}$. 

Suppose that $b\in\C,\ 0<r_1<r_2<|b|$ and suppose that $f$ extends holomorphically 
from each circle $b\D (a,\rho)\subset A(b,r_1,r_2)$ which
surrounds $b$.  Notice that no such 
$b\D (a,\rho)$ surrounds the origin. By Theorem 2.1 the 
function  $F(z,\overline z)=f(z)$ 
has a bounded continuous extension $F_1$ from $\tilde A((b,r_1,r_2)$ to 
$\Omega (A(b,r_1,r_2))\cup b\Omega (A(b,r_1,r_2))$  which is holomorphic on 
$\Omega (A(b,r_1,r_2))$. 

Let $r=(r_1+r_2)/2$ and consider the circle $b\D (b,r)$ and associated varieties 
$\Lambda _{b,r}$ and $V_{b,r}$. Note that $\Omega (A(b,r_1,r_2))$ is 
an open neighbourhood of
$\Lambda _{b,r}$.  Recall that $\Omega (A(b,r_1,r_2))$ is the disjoint 
union of $\Lambda _{a,r}$ such that 
$b\D (a,r)\subset \hbox{Int}A(b,r_1,r_2)$ surrounds $b$.

Write  $V_j = V_j(b,r),\ 1\leq j\leq 4$ and let 
$$
\lambda =\{ (z,\overline{b}+r^2/(z-b))\colon\ z\in b\D(b,r)\cap
\overline \D (0,\sqrt{|b|^2-r^2})\}.
$$
Note that $\lambda $ is an arc which is a part of $b\Lambda _{b,r}\subset b\Omega$. Clearly 
$w\in b\D (b,r)\cap\overline \D (0,\sqrt{|b|^2-r^2})$ if 
and only if $(w,\overline w)\in\lambda$. Each such $w$ is contained in two circles, 
tangent from outside to each other at $w$, one contained in 
$\hbox{Int}A(0,R_1,R_2)$  for some $R_1, R_2,\ 0<R_1<R_2<\infty $, and surrounding the 
origin, and the other contained in $\hbox{Int}A(b,r_1,r_2)$ and surrounding $b$. 
Proposition 2.4 
implies that there are a neighbourhood $U\subset \Sigma$ of $(w,\overline w)$, an open 
convex cone $K\subset i\Sigma $ with vertex at the origin, and an $\eta>0$
such that if
$K_\eta =\{Z\in K,\ |Z|<\eta \}$ then $U+ K_\eta \subset\Omega (A(0,R_1,R_2))$ 
and $U-K_\eta \subset\Omega (A(b,r_1,r_2))$. By 
the edge of the wedge theorem 
it follows that $F(z,\overline z)= f(z)$ has a holomorphic extension $\Phi$ into a small 
open ball $B\subset\C^2$ centered at $(w,\overline w)$. Provided that $B$ is small enough
Proposition 3.1 implies that  
$\Phi\equiv F$ on $B\cap\Omega $ and $\Phi\equiv F_1$ on $B\cap \Omega (A(b,r_1,r_2))$.
Since we can repeat the process for every 
$(w,\overline w)\in\lambda$ it follows that 
there is an open connected
neighbourhood $\cV$ of $\lambda $ in $\C^2$ 
such that that $\cV\cap \Omega$,\ $\cV\cap\Omega (A(b,r_1,r_2))$ and
$\cV\cap \tilde A(b,r_1,r_2)$ 
are connected, such that $F(z,\overline z)=f(z)$ has a holomorphic 
extension $\Phi$ to $\cV$ which satisfies  
$\Phi\equiv F$ on $\Omega\cap \cV$ and $\Phi\equiv F_1$ on 
$\Omega (A(b,r_1,r_2))\cap \cV$. 

The components $V_3$ and $V_4$ of $V_{b,r}\cap\Omega$ are contained 
in $\Omega $ so $F$ is well defined 
and holomorphic on $V_3$ and $V_4$. The function $F_1$ is 
well defined on components $V_1$ and $V_4$ of $V_{b,r}\setminus b\Omega$ 
which  
together with the arc $\{(z,\overline{b}+r^2/(z-b))
\colon\ z\in b\D (0,\sqrt{|b|^2-r^2}\cap\D (b,r)\}$ form    
 $\Lambda _{b,r}$. We first show that 
on $V_4$, where both $F$ and $F_1$ are defined, 
these two functions coincide. To see this, 
choose $w\in b\D (b,r)\setminus\overline\D (0,\sqrt{|b|^2-r^2})$. 
There is a disc $\D (c,R)$ which contains the origin such that $\D(b,r)\subset
\D (c,R)$ and such that $b\D (c,R)$ is tangent to $b\D (b,r)$ at $w$. 
 Proposition 2.4 implies that   
there are an open  
neighbourhood $U\subset\Sigma $ of 
$(w,\overline{w})$, an open convex cone 
$K\subset i\Sigma$ with vertex at the origin and 
an $\eta >0$ such that if $K_\eta =\{Z\in K,\ |Z|<\eta \}$ then 
$U+K_\eta \subset\Omega (A(b,r_1,r_2))\cap\Omega (A(0,R_1,R_2))$ for some 
$R_1, R_2, 0<R_1<R_2<\infty $ and such that $V_4$ meets $U+K_\eta$. This 
 implies that 
$F\equiv F_1$ on $U+K_\eta$ \ since their boundary values 
$F(z,\overline z)=f(z)=F_1(z,\overline z)\ ((z,\overline z)\in U)$, are the same. 
Since $V_4$ meets $U+K_\eta $ it follows that $F\equiv F_1$ on $V_4$.

The arc $\lambda $ is the intersection of the boundaries of $V_1$ and $V_3$ in $V_{b,r}$. 
We show that $F_1|V_1$ is
the analytic continuation of $F|V_3$  in $V_{b,r}$ across
$\hbox{Int}\lambda$. To see this, recall that there are an open 
neighbourhood $\cV\subset \C^2$ of $\lambda$
and a holomorphic function $\Phi  $ on $\cV$ such that  
$\Phi \equiv F$ on $\cV\cap\Omega $ and 
$\Phi\equiv F_1$ on $\cV\cap\Omega (A(b,r_1,r_2))$, so there is 
a single holomorphic function $\Psi = \Phi | V_{b,r}\cap\cV$ on $V_{b,r}\cap\cV$ such that
$\Psi\equiv F$ on $V_3\cap \cV$ and 
$\Psi\equiv F_1$ on $V_1\cap\cV$.  

Thus we showed that $F$ extends holomorphically into a neighbourhood $\cV$ of $\lambda $ in $\C^2$, that 
$F|V_3\cup V_4$ extends holomorphically along $V_{b,r}$ into a neighbourhood of 
$\overline {V_1}$ in $V_{b,r}$ and that $F|V_4\equiv F_1|V_4$. 

We now use the preceding reasoning further to show that $F$ extends holomorphically 
to a neighbourhood of $\overline{V_1}$ in $\C^2$. To see this, we choose a small $\eta >0$ 
and repeat the process above with $V_{a,r},\ a\in \D (b,r),$ in place of $V_{b,r}$. Note that the union 
$\cW$ of all $\Lambda_{a,r},\ a\in\D (b,r)$, is an open neighbourhood of $\Lambda_{b,r}$ which, 
provided that $\eta$ is small enough, is contained in $\Omega(A(b,r_1,r_2))$ and so $F_1$ is 
holomorphic on $\cW$. Note that $\cV\cup\cW$ is a neighbourhood of $\overline{V_1}$ in $\C^2$. 
Repeating the process above for $a\in\D (b,\eta)$ in place of $b$ we see that 
$F|V_3(a,r)\cup V_4 (a,r)$ extends holomorphically along $V_{a,r}$ into a neighbourhood of 
$\overline{V_1(a,r)}$ in $V_{a,r}$. However, in $\cW$ all these extensions coincide with $F_1$ 
so the function $\Psi $ on $\cV\cup\cW$ defined as $\Psi|\cV = F|\cV$ and $\Psi|\cW = F|\cW$ 
then $\Psi $ is holomorphic on $\cV\cap\cW$. Thus, $F$ extends holomorphically into 
$\cV\cap\cW$, a neighbourhood of $\overline{V_1}$ in $\C^2$. 

We will now apply the continuity principle. Recall that $F$ extends holomorphically 
into a neighbourhood 
$P$ of $\overline{V_1}$ in $\C^2$. Now, 
$$
V_1 = \bigl\{ \bigl(z,\overline b+r^2/(z-b)\bigr)\colon\ z\in D_1(b,r)\bigr\}
$$
is an embedded analytic disc whose boundary
$$
bV_1 = \bigl\{ \bigl(z,\overline b+r^2/(z-b)\bigr)\colon\ z\in bD_1(b,r)\bigr\}
$$
is contained in $b\Omega $. For each $t,\ 0\leq t\leq 1$, let 
$$
V_{1,t} = \bigl\{ \bigl(tz,\overline b+r^2/(z-b)\bigr)\colon\ z\in D_1(b,r)\bigr\}.
$$
Then $V_{1,t},\ 0\leq t\leq 1$, is a continuous family of embedded analytic discs, 
$V_{1,1}=V_1$, whose boundaries
$$
bV_{1,t} = \bigl\{ \bigl(tz,\overline b+r^2/(z-b)\bigr)\colon\ z\in bD_1(b,r)\bigr\}
$$
are contained in $\Omega\cup P$\ (in fact, for $0\leq t <1$ they are contained in $\Omega $). 
By the continuity principle it follows that $F$ extends holomorphically into a neighbourhood 
of 
$$
\overline{V_{1,0}} = \bigl\{ \bigl(0,\overline b+r^2/(z-b)\bigr)\colon\ z\in 
\overline{D_1(b,r)}\bigr\}
$$
in $\C^2$. It is easy to see that $V_{1,0}$ contains the origin. Consequently $F$ extends 
holomorphically into a neighbourhood of the origin in $\C^2$ and so $f$ extends across 
the origin in $\C$ as a function which is real analytic in a neighbourhood of the origin. Since $f$ 
extends holomorphically from every circle surrounding the origin it follows from [G2] that $f$ 
is holomorphic on $\C$. This completes the proof.
\vskip 2mm
\noindent \bf Remark\ \rm In the last step of the proof above we may, instead of [G2],
use the Liouville theorem as follows: From (2.2) 
it follows that $F$ is bounded on $
\Delta (0,\delta )\times\C $ for some $\delta >0$. By the 
Liouville theorem the function 
$\z\mapsto F(z,\z )$ is constant for each 
$z\in\Delta (0,\delta )$ so $F$ does not 
depend on $w$ on   $\Delta (0,\delta )\times\C $. 
It follows that $F$ is a function of $z$ only so  
$f(z)=F(z,\overline z)\ \ (z\in\C
\setminus\{0\})$ is a restriction of an entire 
function to $\C\setminus\{0\}$.
\vskip 3mm
\bf 7.\ Examples\rm 
\vskip 2mm
By Theorem 1.1 the family of all circles that 
surround the origin is a maximal \it open \rm family of 
circles that is not a test family for holomorphy. Even its closure, that is 
the family of all circles that 
either surround the origin or pass through the origin is not a maximal family 
that is not a test 
family for holomorphy. To see this, let $a\in\C,\ \rho >0, |a|>\rho$, and let 
$$
g(z)=
\left\{\eqalign{&0\hbox{\ \ if\ }z=0\cr
                &(z^2/\overline z)[(z-a)
				(\overline z-\overline a)-\rho^2]\hbox{\ \ if\ }z\not=0.\cr}\right.
				$$
				The function $g$ vanishes identically on $b\D (a,\rho)$ 
				and hence extends holomorphically from 
$b\D (a,\rho)$. Since $(z^2/\overline z)[(z-a)(\overline z-\overline a)-\rho^2]= 
z^3-az^2-\overline a (z^3/\overline z) 
+a\overline a (z^2/\overline z)-\rho^2 (z^2/\overline z)$ is a polynomial
in $z$ and $1/\overline z$ it follows that $g$ 
extends holomorphically from every circle that surrounds the origin. 
Since $g$ is continuous on $\C$ it
extends holomorphically also from every circle that passes through
the origin. This shows that if 
$|a_i|>\rho_i>0,\ 1\leq i\leq n$, then the function
$$
f(z) =
\left\{\eqalign
{&0\hbox{\ \ if\ }z=0\cr
&(z^2/z)^n\Pi_{j=1}^n[(z-a_j)(\overline z-\overline a_j)-\rho _j^2]\hbox{\ \ if\ }z\not=0
\cr}\right.
$$ 
is continuous on $\C$, extends holomorphically 
from all circles that surround the origin, from all 
circles that pass through the origin, and 
from all circles $b\D (a_i,\rho_i),\ 1\leq i\leq n$ yet $f$ is not holomorphic.

In our next example, let $g$ be a function from the disc algebra and define 
$$
f(z)=g(z/\overline z)\ \ (z\in\C\setminus\{0\}).
\eqno (7.1)
$$
Suppose that $b\D (a,\rho)$ surrounds the origin. Then $|a|<\rho$ and for $|\z|=1$ we have 
$$
f(a+\z\rho)= g\biggl({{a+\z\rho}\over{\overline a+\overline\z\rho}}\biggr)=
g\bigg(\z{{\z+a/\rho}\over{1+(\overline a/\rho)\z}}\biggr)
$$
which shows that the function $\z\mapsto f(a+\z\rho)\ \ (\z\in b\D )$ 
extends to a function from the disc algebra. Thus, $f$ extends holomorphically 
from every circle surrounding the origin. Since the boundary values of the functions 
from the disc algebra can be highly non-smooth this example shows that 
a highly nonsmooth function on 
$\C\setminus \{0\}$ can be holomorphically extendible 
from every circle surrounding the origin.  
\vskip 3mm
\bf 8.\ Analyticity on circles for functions constant on lines\rm
\vskip 2mm
In the second example in Section 7 the function $f$ is 
constant on each line passing through the origin,
 that 
 is, 
 $$
 f(tz)=f(z) \ \ (z\in\C\setminus\{ 0\},\ t\in\R\setminus\{ 0\}).
 \eqno (8.1)
 $$
 In this section we look more closely at such functions.
 \vskip 2mm
 \noindent \bf Theorem 8.1\ \it Suppose that $f$ is a
 continuous function on $\C\setminus\{ 0\}$ 
 that is a constant on each line passing through the origin, 
 that is, $f$ satisfies (8.1). 
 If $f$ extends holomorphically from one circle surrounding 
 the origin then it extends holomorphically 
 from every circle surrounding the origin. This happens if 
 and only if there is a function $g$ 
 from the disc algebra such that $f(z)=g(z/\overline z)\ 
 (z\in\C\setminus\{ 0\})$. \rm 
 \vskip 1mm
 \noindent\bf Proof.\ \rm Suppose that $f$ is a continuous 
 function on $\C\setminus\{0\}$ that satisfies (8.1). 
 Then there is a continuous function $g$ on $b\D $ such that
 $$
 f(z)=g(z/\overline z)\ \ (z\in\C\setminus\{0 \}.
 \eqno (8.2)
 $$
 Assume that $f$ extends holomorphically from a circle $b\D (a,\rho)$ 
 that surrounds the origin. By (8.1) we may assume 
 that $\rho = 1$. If $a=0$ then the function 
 $\z\mapsto g(\z^2)\ (\z\in b\D)$ extends to a 
 function in the disc algebra which implies that $g$ 
 extends to a function from the disc algebra. 
 Suppose that $a\not= 0$. Composing $f$ with a rotation 
 if necessary we may assume that 
 $0<a<1$. By our assumption there is a function $h$ from 
 the disc algebra such that 
 $$
 h(\z )= g\biggl({{a+\z }\over{a+\overline\z}}\biggr)=
 g\biggl(\z{{a+\z}\over{1+\z a}}\biggr)\ \ (\z\in b\D).
 $$
 Clearly $\xi\mapsto h\bigl((\xi-t)/(1-t\xi)\bigr)\ 
 \ (\xi\in b\D)$ extends to a function from 
 the disc algebra for every $t,\ 0\leq t <1$. Put 
 $t= (1-\sqrt{1-a^2})/a$. Then $0<t<1$ and 
 $$
 {{a+{{\xi -t}\over{1-\xi t}}}\over{1+a{{\xi -t}\over{1-\xi t}}}}=
 {{\xi+t}\over{1+t\xi}}\ \ \ (\xi\in b\D)
 $$
 which implies that
 $$
 \xi\mapsto g\biggl({{\xi -t}\over{1-\xi t}}{{\xi +t}\over{1+\xi t}}\biggr)
 = g\biggl({{\xi^2-t^2}\over{1-t^2\xi^2}}\biggr)\ \ (\xi\in b\D)
 $$
 extends to a function from the disc algebra which implies that 
 $\xi\mapsto  g((\xi-t^2)/(1-t^2\xi))\ \ (\xi\in b\D )$ 
 extends to a function from the disc algebra. Consequently $g$ extends to 
 a function in the disc algebra and so $f$ is of the form (7.1). By the discussion 
 following (7.1) it follows that the function $f$ extends holomorphically 
 from every circle that surrounds the origin.
 \vskip 3mm
 \bf 9.\ More examples \rm
 \vskip 2mm
 \noindent \bf Example 9.1\ \ \rm Let $0<a<1$ and let $\Phi (\z )=(a+\z )/|a+\z |\ \ 
 (\z\in b\D )$. Then $\Phi\colon\ b\D\rightarrow b\D$ is diffeomorphism. Define
 $$
 f(z) = \Phi^{-1}(z/|z|)\ \ (z\in\C\setminus \{ 0\}).
 \eqno (9.1)
 $$
 The function $f$ is continuous on $\C\setminus \{ 0\}$ and is constant 
 on each ray emanating from the origin, that is,
 
 $$
 f(tz)=f(z)\ \ (z\in\C\setminus \{ 0\},\ t>0).
 \eqno (9.2)
 $$
 If an $f$ satisfying (9.2) extends holomorphically from a 
 circle $b\D (a,\rho)$ 
 that surrounds the origin then it extends holomorphically
 from $b\D (ta,t\rho)$ for every $t>0$. 
 So, when studying holomorphic extendibility from circles 
 $b\D (a,\rho)$ we may, with no loss of generality, 
 assume that $\rho =1$. 
 
 Let $f$ be as in (9.1). Since 
 $$
 f(a+\z ) = \Phi ^{-1} ((a+\z )/|a+\z |) = \z\ \ (\z\in b\D)
 $$
 it follows that $f$ extends holomorphically from $b\D (a,1)$. By (9.2) 
 $$
 {{a+\z}\over{|a+\z|}} =
 \sqrt{{{(a+\z)^2}\over{(a+\z)(a+\overline\z )}}} =
 \sqrt{{a+\z}\over{a+\overline\z}} =
 \sqrt{\z {{a+\z}\over{1+a\z }}}\ \ (\z\in b\D )
 $$
 it follows that
 $$
 \z\mapsto f\biggl(\sqrt{\z{{a+\z}\over{1+a\z }}}\biggr)
 $$
 extends to a function from the disc algebra which happens if and only if
 
$$
 \z\mapsto f\biggl(\sqrt{M(\z ){{a+M(\z )}\over{1+aM(\z )}}}\biggr)
 $$
 extends to a function from the disc algebra for an automorphism $M$ of $\D$. 
 In particular, 
 if $M(\z )= (\z-a)/(1-a\z)$ it follows that 
 $$
 \z\mapsto f\biggl(\sqrt{\z {{\z-a}\over{1-a\z }}}\biggr)\ \ (z\in b\D )
 $$
 extends to a function from the disc algebra which is 
 equivalent to the fact that $f$ extends 
 holomorphically from $b\D (-a,1)$. It will follow from 
 Theorem 10.1 that these two circles 
 are the only circles  of radius $1$ from which $f$ 
 extends holomorphically.
 \vskip 2mm
 \noindent\bf Example 9.2\ \ \rm Let $g$ be a function 
 in the disc algebra which is  not an 
 even function. Let  $f(z) = g(z/|z|)\ \ (z\in\C\setminus \{ 0\})$. 
 Then $f$ is continuous on 
 $\C\setminus\{ 0\}$ and extends holomorphically from $b\D (0,1)$. 
 It will follow from 
 Theorem 10.1 below that $b\D (0,1)$ 
 is the only circle of radius $1$ from which $f$ extends holomorphically. 
 \vskip 3mm
 \bf 10.\ Analyticity on circles for functions constant on rays \rm
 \vskip 2mm
 In both examples in Section 9 the function $f$ satisfies (9.2), that is, 
 $f$ is constant on each ray emanating from the origin. In this section we 
 look more closely at such functions.
 
 Suppose that a continuous function $f$ on 
 $\C\setminus\{ 0\}$ satisfies (9.2). Assume that 
 $0\leq d<1$ and that $f$ extends holomorphically 
 from  $b\D (de^{i\alpha}, 1)$ for some $\alpha\in\R$. 
 This means that $\z\mapsto f(e^{i\alpha}(d+\z ))\ \ (\z\in b\D)$ 
 extends to a function in the disc algebra. 
 Since 
 $$
 e^{i\alpha}(d+\z )/|d+\z |=e^{i\alpha}\sqrt{\z (d+\z )/(1+d\z ) }\ \ (\z \in b\D)
 $$
 this happens if and only if
 $$
 \left.\eqalign{
 f\bigl(e^{i\alpha}\sqrt{\z (d+\z )/(1+d\z ) }\bigr) = q(\z )\ \ (\z\in b\D )
 \cr
 \hbox{where \ } q\hbox{\ belongs to the disc algebra.}
 \cr}
 \right\}
 \eqno(10.1)
 $$  
 In the case when $d=0$ this implies that $\z\mapsto f(\z )\ (\z\in b\D)$ 
 extends to a function 
 from the disc algebra. Suppose that $d\not= 0$. Put
 $$
 \z = {{\xi - t}\over{1-t\xi}}\hbox{\ \ where\ \ } t={{1-\sqrt{1-d^2}}\over{d}}
 $$
 to get
 $$
 \z{{d+\z }\over{1+d\z}}={{\xi^2-t^2}\over{1-t^2\xi^2}}
 $$
 so that (10.1) is equivalent to 
 $$
 \left.\eqalign{
 f\biggl(e^{i\alpha}\sqrt{{{\xi^2-t^2}\over{1-t^2\xi^2}} }\biggr) 
 = G(\xi )\ \ (\xi\in b\D )
 \cr
 \hbox{where \ } G\hbox{\ belongs to the disc algebra.}
 \cr}
 \right\}
 \eqno(10.2)
 $$
 In fact, $G(\xi )=q((\xi-t)/(1-t\xi))\ (\xi\in b\D )$. Putting $Z=e^{i\alpha}\sqrt{
 (\xi^2-t^2)/(1-t^2\xi^2)}$ we get
 $$
 \xi^2={{(e^{-i\alpha}Z)^2+t^2}\over{1+t^2(e^{-i\alpha}Z)^2}}
 $$
 which implies that (10.2) is equivalent to 
 $$
 \left.\eqalign{
 f(Z) = G\biggl(
 e^{-i\alpha}
 \sqrt{
 {Z^2+e^{2i\alpha} t^2}
 \over
 {1+e^{-2i\alpha}t^2Z^2} 
 }
 \biggr) \ \ (Z\in b\D )
 \cr
 \hbox{where \ } G\hbox{\ belongs to the disc algebra.}
 \cr}
 \right\}
 \eqno(10.3)
 $$
 \vskip 2mm
 \noindent{\bf Theorem 10.1}\ \ \it Let $f$ be a continuous function on 
 $\C\setminus\{0\}$ which satisfies (9.2), that is, $f$ is constant 
 on each ray emanating from the origin. Suppose that $f$ extends holomorphically from 
 $b\D (a,1)$ and $b\D (b,1)$ where $a,b\in\D ,\ b\not=a, b\not=-a$. 
 Then there is a function $g$ in the disc algebra such that 
 $$
 f(z)=g(z/|z|)\ \ (z\in\C\setminus\{ 0\}).
 \eqno (10.4)
 $$
 Consequently, $f$ satisfies (8.1), that is, $f$ is 
 constant on each line passing through the origin
 and extends holomorphically from each circle surrounding the origin. \rm
 \vskip 1mm
 \noindent \bf Proof.\ \ \rm Suppose that $f$ extends holomorphically 
 from $b\D (d_1e^{i\alpha_1},1)$ and $b\D (d_2e^{i\alpha_2},1)$  where $0\leq d_i<1\ (i=1,2), 
 \ d_2e^{i\alpha_2}\not=d_1e^{i\alpha_1},\ d_2e^{i\alpha_2}\not=-d_1e^{i\alpha_1}$. 
 It is enough to prove that $f$ is an even function for then the rest 
 follows from Theorem 8.1. 
 
 Let $t_i=0$ if $d_i=0$ and $t_i = (1-\sqrt{1-d_i^2})/d_i$ if $d_i\not= 0,\ i=1,2$. Write
 $A_i=e^{2i\alpha_i}t_i^2,\ i=1,2$. By the discussion 
 preceding Theorem 10.1 there are functions $G_1, G_2$ in the disc algebra such that 
 $$
 f(Z)=G_i\bigl(e^{-i\alpha_i}\sqrt{{{Z^2+A_i}\over{1+\overline{A_i}Z^2}}}\bigr)
 \ \ (Z\in b\D,\ i=1,2)
 \eqno (10.5)
 $$
 Write $W^2= (Z^2+A_1)/(1+\overline{A_1}Z^2)$ so that $Z^2=(W^2-A_1)/(1-\overline{A_1}W^2)\ \ 
 (W\in b\D)$ and $(Z^2+A_2)/(1+\overline{A_2}Z^2)=(W^2+C)/1+\overline C W^2)$ 
 where $C= (A_2-A_1)/(1-A_1 A_2)$.    
 Now (10.5) implies that 
 $$
 G_1\biggl(e^{-i\alpha _1}\sqrt{{Z^2+A_1}\over{1+\overline{A_1}Z^2}}\biggr)
 =
 G_2\biggl(e^{-i\alpha _2}\sqrt{{Z^2+A_2}\over{1+\overline{A_2}Z^2}}\biggr)
 \ \ (Z\in b\D)
 $$
 which implies that 
 $$
 G_1(e^{-i\alpha _1}W)= G_2\bigl(e^{-i\alpha_2}\sqrt{{{W^2+C}\over{1+\overline CW^2}}}
 \bigr)\ \ (W\in b\D).
 \eqno (10.6)
 $$
 Since both $G_1 $ and $G_2$ belong to the disc algebra it follows that the 
 relation (10.6) continues holomorphically inside $\D $, 
 so (10.6) implies that either $C=0$ or $G_2$ is an even 
 function. Assume that $C=0$. By (10.6) it follows that $A_1=A_2$ and 
 $e^{i\alpha_1}=e^{i\alpha _2}$. It follows that 
 $d_2e^{i\alpha_2}=\pm d_1e^{i\alpha_1}$ 
 which is impossible by the assumption. 
 Thus, $G_2$ is an even function and consequently, by (10.5), 
 $f$ is an even function.  
 This completes the proof. 
 \vskip 2mm
 \noindent\bf Remark\rm\ \ Note that Theorem 10.1 implies that in Example 9.1 
 the circles $b\D (a,1)$ and 
 $b\D (-a,1)$ are the only circles of radius one 
 from which $f$ extends holomorphically. Similarly, 
 in Example 9.2, \ $b\D $ is the only circle of radius 
 one from which $f$ extends holomorphically.
 \vskip10mm
 \noindent\bf Acknowledgement\ \rm The fundamental observation in [AG] that the 
 holomorphic extensions of rational 
 functions from circles are naturally defined on varieties $\Lambda _{a,\rho}$, objects 
 in $\C^2$, is due to 
 Mark Agranovsky. The author is grateful to him for many stimulating 
 discussions about holomorphic extensions from circles.
 
 This work was 
 supported in part by the Ministry of Education, 
 Science and Sport of the Republic of Slovenia. 
\vfill
\eject
\centerline{REFERENCES}
\vskip 5mm
\noindent[AG]\ M.\ Agranovsky, J.\ Globevnik:\ Analyticity on 
circles for rational and real analytic functions of two 
real variables.

\noindent To appear in J.\ d'Analyse Math.
\vskip 2mm
\noindent [G1]\ J.\ Globevnik:\ Testing analyticity on rotation invariant families 
of curves.

\noindent Trans.\ Amer.\ Math.\ Soc.\ 306 (1988) 401-410
\vskip 2mm
\noindent [G2]\ J.\ Globevnik:\ Integrals over circles passing 
through the origin and a characterization of analytic functions.

\noindent J.\ d'Analyse Math. 52 (1989) 199-209
\vskip 2mm
\noindent [G3]\ J.\ Globevnik:\ Holomorphic extensions from open families of circles.

\noindent Trans.\ Amer.\ Math.\ Soc.\ 355 (2003) 1921-1931
\vskip 2mm
\noindent [T]\ A.\ Tumanov:\ A Morera type theorem in the strip.

\noindent To appear in Math.\ Res.\ Lett.\ 
\vskip 20mm
\noindent Institute of Mathematics, Physics and Mechanics

\noindent University of Ljubljana, Ljubljana, Slovenia

\noindent josip.globevnik@fmf.uni-lj.si

\bye